\documentclass[review]{elsarticle}
\UseRawInputEncoding
\usepackage{amsfonts,amssymb,amsmath,amsbsy,amsthm,amscd}
\usepackage[english]{babel}
\usepackage{graphicx}

\newtheorem{theorem}{Theorem}[section]

\newtheorem{lemma}[theorem]{Lemma}
\newtheorem{corollary}{Corollary}
\newtheorem{proposition}{Proposition}

\renewcommand{\le}{\leqslant}
\renewcommand{\ge}{\geqslant}

\begin{document}

\begin{frontmatter}

\title{Localization of the formation of singularities in multidimensional compressible Euler equations }


\author[1]{Olga S. Rozanova}
\ead{rozanova@mech.math.msu.su} 

\address[1]{Department of Mechanics and Mathematics, Moscow
State University, Moscow 119991 Russia}


\begin{abstract}
We consider the Cauchy problem with smooth data for compressible Euler equations in many dimensions and concentrate on two cases: solutions with finite mass and energy and solutions corresponding to a compact perturbation of a nontrivial stationary state.
We prove the blowup results using the characteristics of the propagation of the solution in space and find upper and lower bounds for the density of a smooth solution in a given region of space in terms of the initial data.
 To solve the problems, we introduce a special family of integral functionals and study their temporal dynamics.
\end{abstract}

\def\sign{\mathop{\rm sgn}\nolimits}






\end{frontmatter}


\section{Introduction}
We consider a system for density $\varrho({\bf x},t)$, pressure
$p({\bf x},t)$, velocity ${\bf V}({\bf x},t)$ and entropy $S({\bf
x},t)$:
\begin{eqnarray}
 \partial_t (\varrho{\bf V}) +   {\rm div}_{\bf x} ({\varrho\bf
V}\otimes{\bf V}) +
 \nabla_{\bf x} p  &=&  0, \label{e1}\\
\partial_t \varrho +  {\rm div}_{\bf
x} (\varrho {\bf V}) &=& 0, \label{e2} \\
\partial_t S + {\bf V}\cdot\nabla_{\bf x} S&=&0,\label{e3}\\
p=e^S \,\varrho^\gamma. \label{e4}
\end{eqnarray}


Here ${\bf x}\in{\mathbb R}^n,$ $t\ge 0$, $\gamma>1$
is the heat ratio.

The initial data are the following:
\begin{eqnarray}\label{cd}
({\bf V},\varrho, p)|_{t=0}=({\bf V}_0,\varrho_0\ge 0, p_0\ge 0)({\bf
x})\label{CD}\in C^1({\mathbb R}^n).
\end{eqnarray}

\medskip
It is well known that system \eqref{e1} -- \eqref{e4} can be written in a symmetric form, therefore the Cauchy problem \eqref{e1} -- \eqref{cd} has locally in time a solution as smooth as initial data (e.g.\cite{Kato}).
At the same time, the determination of the class of initial data leading to the formation of a singularity is a very interesting and important problem from both theoretical and practical points of view.
Since Sideris' seminal paper, \cite {Sideris}, there have been many papers in which integral moments of mass are used to find sufficient conditions for the blowup of initially smooth solutions to the Cauchy problem in the multidimensional case (for an overview, see \cite {Chen}, \cite {Xin}, \cite {Rozanova_JMS}). Almost all of them use properties of the integral $G(t)=\frac12 \int\limits_{{\mathbb R}^n}\rho |{\bf x}|^2 \, dx$ (see \cite{Steven}, \cite{Rozanova_DU} for other variants, the case of $G(t)$ is our particular case, see Corollary \ref{cor1}). Moreover, the previous methods do not give an answer about the position of the emerging singularity in space. The present work seems to be a first attempt to use compactly supported moment function and describe the properties of the smooth solution of the full compressible Euler equations inside a fixed space domain. In addition, we use this new kind of moment function to find sufficient conditions for the behavior of the formation of a solution singularity and to determine the necessary behavior of a globally smooth solution in time as $ | {\bf x} | \to \infty $ (see \cite {Steven1} for another kind of compactly supported moment functions applied to a simpler situation).

We consider two cases.

I.  Solutions with finite total mass
$$m=\int\limits_{{\mathbb R}^n}\rho\,dx={\rm const}>0$$ and finite total energy
$$\mathcal E=E_k(t)+E_p(t)=\frac{1}{2}\int\limits_{{\mathbb R}^n}\rho|{\bf V}|^2
\,dx+\frac{1}{\gamma-1}\int\limits_{{\mathbb R}^n} p\,dx={\rm const}>0.$$

II. Perturbations of the stationary state with initial data
\begin{eqnarray}\label{CD2}
\quad (\varrho_0, p_0)({\bf
x})>0,\quad ({\bf V}_0,\varrho_0, p_0)({\bf x})=({\bf 0},
\bar\varrho, \bar
p) \,\, \mbox{for}\,\, {\bf x}\notin {\mathcal B}_{{\mathcal R}_0}, \\
\bar\varrho, \bar p ={\rm const}>0,\, {\mathcal B}_{{\mathcal R}_0}={\mathcal B}_{{\mathcal R}(0)}=\{{\bf x} \big||{\bf x}|<{\mathcal R}_0,\,
{\mathcal R}_0={\rm const}>0\}.\nonumber
\end{eqnarray}

In the latter case, we consider
\begin{eqnarray*}
m(t) = \int_{{\mathcal B}_{{\mathcal R}(t)}} (\varrho-  \bar\varrho) \, dx,\quad
e(t) = E_k(t)+\frac{1}{\gamma-1} \int_{{\mathcal B}_{{\mathcal R}(t)}} (p- \bar p) \, dx,
 \end{eqnarray*}
where ${\mathcal B}_{{\mathcal R}(t)}=\{{\bf x}\big| |{\bf x}|<{\mathcal R}(t)={\mathcal R}_0+\sigma
t\}$, $\sigma=\sqrt{p_\varrho}|_{\varrho=\bar\varrho}$ is the sound
speed (speed of propagation of perturbations). We also denote
$${E_p}_{{\mathcal B}_{\mathcal R}(t)}=\frac{1}{\gamma-1} \int_{{\mathcal B}_{{\mathcal R}(t)}} p \, dx,\quad {\mathcal E}_{{\mathcal B}_{{\mathcal R}(t)}}=E_k(t)+{E_p}_{{\mathcal B}_{{\mathcal R}(t)}}.$$

 Recall ( \cite{Sideris}) that for $C^1$ - smooth solutions of \eqref{e1} -- \eqref{CD} the support of perturbation is contained in  ${\mathcal B}_{{\mathcal R}(t)}$ and
 \begin{eqnarray}\label{dm}
  m'(t)=0,\quad
e'(t)=0.\nonumber
 \end{eqnarray}

\bigskip

We consider the functional
\begin{equation}
G_\phi(t)=\int\limits_{{\mathbb R}^n}\rho \phi(|{\bf
x}|)\,dx.\label{Gphi} \end{equation}

Let us introduce $\phi(|{\bf x}|) \in C^1([0,\infty))$ as follows:
\begin{equation}\label{phi}
 \phi(|{\bf
x}|)=\left\{
\begin{array}{l}1-B |{\bf x}|^2,\quad |{\bf x}|\in [0,R-\frac{R}{k})
\\
C(|{\bf x}|-R)^2,\quad |{\bf x}|\in
(R-\frac{R}{k},R]\qquad\qquad\\0,\quad |{\bf x}|\in [R,+\infty)
\end{array}\right.
\end{equation}
Here $R>0,$ $k>n$,
$B=\frac{k}{R^2(k-1)},\,$ $C=B(k-1),$ the support of $\phi(|{\bf x}|)$
is ${\mathcal B}_R=\{x\Bigl||{\bf x}|\le R\}.$

We denote
$$Q_\phi(t)=\int\limits_{{\mathbb R}^n}(\rho-\bar \rho) \phi(|{\bf
x}|)\,dx,\quad K=\int\limits_{{\mathbb R}^n} \phi(|{\bf
x}|)\,dx,$$
$$\delta=\min\{2, n(\gamma-1)\}>0,$$
$$A_1:=(k-n)\exp\{\inf\limits_{{\mathcal B}_R}
S_0\}\left(\int\limits_{{\mathcal B}_R}\phi^{\frac{\gamma}{\gamma-1}}\,dx\right)^{1-\gamma},\quad
A_2:=\max\{2,(\gamma-1)k\}, $$
 $\omega_n$  is the volume of $n$ - dimensional unit ball.

We are going to prove the following theorems.

\medskip

{\bf Case I.}

\begin{theorem}\label{T1}  The solution of problem \eqref{e1}--\eqref{e4}, \eqref{cd}  with a finite total energy $\mathcal E$ cannot keep initial smoothness for all $t>0$, if
\begin{equation}\label{condT0}
\lim_{R\to\infty}  {R^n} \sup\limits_{{|x|=R}}\left[\frac{\delta}{R}\int_0^t\left(\frac12 \rho|{\bf V}|^2+ \frac{\gamma}{\gamma-1} p\right)|{\bf V}|\,dt +(\rho|{\bf V}|^2+n p)\right]\le {\delta_1 {\mathcal E}},
\end{equation}
 with any constant $\delta_1,$ such that $\displaystyle 0\le \delta_1<\frac{\delta}{n\omega_n}$.
\end{theorem}

\medskip

\begin{theorem}\label{T2}
Let $z(t)$ be the solution of problem
\begin{equation}
 z''(t)= 2B\left(A_1 z^\gamma(t)-A_2\mathcal
E\right), \quad z(0)=G_\phi(0), \, z'(0)=G'_\phi(0),\label{(1.6.5e)}
\end{equation}
with a fixed  $R$ and $k$ (introduced in  \eqref{phi}) such that $z(t)\ge 0$ for $t\in [0, T]$. For
 all $t\in [0, T]$ such that the solution to the Cauchy problem  \eqref{e1}--\eqref{e4}, \eqref{cd} keeps
initial smoothness inside ${\mathcal B}_R$, the following estimates hold:
 \begin{eqnarray}
\sup\limits_{x\in B_R} \rho(t,x)&\ge& \frac{z_-(t)}{K},\label{T21}\\
\inf\limits_{x\in B_R} \rho(t,x)&\le &\frac{z_+}{K},\label{T22}
\end{eqnarray}
where $z_-(t)=z(t)$, $z_+=\min\{m, \left((\gamma-1)(k-n)\frac{\mathcal E}{A_1}\right)^\frac{1}{\gamma}   \}$.
If  \eqref{T21} fails for some $t_*\in [0, T]$  or  \eqref{T22} fails for some $t_*>0$, then the solution looses smoothness within $(0,t_*) $. The value of $T$ is given by \eqref{cond}, \eqref{time1}, \eqref{time2}.
\end{theorem}

\medskip
\newpage
{\bf Case II.}

 \medskip

 \begin{theorem}\label{T3}
 Denote $$
 N=\delta e(0)+ \omega_n R^n \bar p\left(\frac{k-1}{k}\right)^n \left(n-\frac{\delta}{\gamma-1}\right).
 $$
 Let constants ${\mathcal R}_0$, $R$ and $k$ be such that ${\mathcal R}_0<\frac{(k-1)R}{k}$ and the
 initial data \eqref{CD2} be such that
 \begin{equation}\label{T3cond}
 0 < T_1 \le  T_2,
    \end{equation}
   with
     $$ T_1=\frac{1}{2 B N}\left(G_\phi'(0)+\sqrt{(G_\phi'(0))^2+4BNG_\phi(0)}\right),\quad
  T_2=\frac{1}{\sigma}\left(\frac{(k-1)R}{k}-  {\mathcal R}_0\right).$$
Then the solution to problem
 \eqref{e1}--\eqref{e4}, \eqref{CD2} loses its smoothness within time $T_1$ inside the ball ${\mathcal B}_R$.
  \end{theorem}
 \medskip
 \begin{theorem}\label{T4}
 Denote $$Q_\phi^+(t)=m(0)-\bar\rho K+\bar\rho \omega_n ({\mathcal R}_0+\sigma t)^n,$$ and let $Q_\phi^-(t)$ be the solution to  problem
 \begin{equation*}
z''(t)= \varkappa^2 z(t)     +P_n(t),\quad z(0)=G_\phi(0)-\bar \rho K, \quad z'(0)=G'_\phi(0),
\end{equation*}
 with $\varkappa^2 = 2\gamma B A_1 (\bar\rho K)^{\gamma-1},$ $P_n(t)= 2 B (A_1 (\bar\rho K)^{\gamma} -A_2 (e(0)+\bar p \,\omega_n \,({\mathcal R}_0+\sigma t)^n)),$ $G_\phi$, $G'_\phi$ are defined in \eqref{Gphi},  \eqref{phi}, \eqref{Gprphi}.
For all $t$ such that the solution to the Cauchy problem  \eqref{e1}--\eqref{e4}, \eqref{CD2} keeps
initial smoothness inside ${\mathcal B}_R$, the following estimates hold:
 \begin{eqnarray}
\sup\limits_{x\in B_R} \rho(t,x)&\ge& \bar\rho +\frac{Q_-(t)}{K},\label{T41}\\
\inf\limits_{x\in B_R} \rho(t,x)&\le &\bar\rho +\frac{Q_+(t)}{K}.\label{T42}
\end{eqnarray}
If  \eqref{T41} or  \eqref{T42} fails for some $t_*>0$, then the solution loses its smoothness within $(0,t_*)$.

\end{theorem}

 \bigskip

The rest of the paper is organized as follows. Section \ref {S2} contains auxiliary results needed for further proof. In Sections \ref {S3} - \ref {S6}, we present the proofs of Theorems \ref {T1} - \ref {T4}, respectively, and give additional comments. Section \ref {S7} discusses the details of the proof and the prospects for using the results obtained.

\section{ Auxiliary results}\label{S2}

\subsection{The main lemma}
First of all, we prove a lemma on the general functional
$$G_\Phi(t)=\int\limits_{{\mathbb
R}^n}\rho(t,x)\Phi(|{\bf x}|)\,dx,$$ considered for all functions
$\Phi(|{\bf x}|)\in C^1[0,+\infty),$ ensuring the convergence of the integral, and such that
$$\lim_{r\to \infty}\int\limits_{S(r)}
 \rho{\bf V}\Phi(|{\bf x}|)\,dS(r)=0,$$
 where $S(r)$ is the $(n-1)$ is the $n$ - dimensional sphere of radius
 $r$.

Let us denote
${\bf
\sigma}=(\sigma_1,...,\sigma_K)$  the vector\,
$\sigma_k=V_i x_j-V_j x_i,
\,i>j,\,i,j=1,...,n,\,k=1,...,K,\,K={\rm C}_n^2.$

\begin{lemma}\label{L1} Provided that all the considered integrals converge, the following equalities hold on the solutions of system \eqref {e1} -
\eqref{e4}:
\begin{equation}\label{Gprphig}
 G'_\Phi(t)=\int\limits_{{\mathbb R}^n}\frac{\Phi'(|{\bf x}|)}{|{\bf x}|}
({\bf V}\cdot{\bf x})\rho\,dx,
\end{equation}

$$G''_\Phi(t)=I_{1,\Phi}(t)+I_{2,\Phi}(t)+I_{3,\Phi}(t)+I_{4,\Phi}(t),$$
where
$$I_{1,\Phi}(t)=\int\limits_{{\mathbb R}^n}\frac{\Phi''(|{\bf x}|)}{|{\bf
x}|^2} |({\bf V}\cdot{\bf x})|^2\rho\,dx,$$
$$I_{2,\Phi}(t)=\int\limits_{{\mathbb R}^n}\frac{\Phi'(|{\bf x}|)}{|{\bf
x}|^3} |{\bf \sigma}|^2\rho\,dx,$$
$$I_{3,\Phi}(t)=\int\limits_{{\mathbb R}^n}\left(\Phi''(|{\bf x}|)+
(n-1)\frac{\Phi'(|{\bf x}|)}{|{\bf x}|} \right) \,p\,dx,$$
$$I_{4,\Phi}(t)=-\lim_{r\to \infty}\int\limits_{S(r)}\Phi'(|{\bf x}|)
 p\,dS(r).$$
\end{lemma}
\proof
The lemma is proved by direct calculation and application of the general Stokes formula. For example, using \eqref{e2}, we obtain
$$G'_\Phi(t)=\int\limits_{{\mathbb R}^n}\rho'_t(t,x)\Phi(|{\bf
x}|)\,dx=-\int\limits_{{\mathbb R}^n}{\rm div}(\rho{\bf
V})\Phi(|{\bf x}|)\,dx=$$$$=\int\limits_{{\mathbb
R}^n}(\nabla\Phi(|{\bf x}|)\cdot {\bf V})\rho\,dx-\lim_{r\to
\infty}\int\limits_{S(r)}
 \rho{\bf V}\Phi(|{\bf x}|)\,dS(r)=\int\limits_{{\mathbb R}^n}\frac{\Phi'(|{\bf
x}|)}{|{\bf x}|} ({\bf V}\cdot{\bf x})\rho\,dx.$$
$\Box$
\medskip

\begin{corollary}\label{cor1} In a particular case $\Phi(|{\bf
x}|)=\frac{|{\bf x}|^2}{2}$,  $|{\bf x}|^2\rho\in L_1({\mathbb R}^n)$,
$$
G'_\Phi(t)=\int\limits_{{\mathbb R}^n}
({\bf V}\cdot{\bf x})\rho\,dx,$$
$$I_{1,\Phi}(t)=\int\limits_{{\mathbb R}^n}\frac{|({\bf V}\cdot{\bf x})|^2
}{|{\bf x}|^2} \rho\,dx,$$
$$I_{2,\Phi}(t)=\int\limits_{{\mathbb R}^n}\frac{|{\bf \sigma}|^2}
{|{\bf x}|^2} \rho\,dx,$$
$$I_{3,\Phi}(t)=n\int\limits_{{\mathbb R}^n}p\,dx=n(\gamma-1)E_p(t),$$
$$I_{4,\Phi}(t)=0.$$
Moreover,
$$I_{1,\Phi}(t)+I_{2,\Phi}(t)=2E_k(t).$$
\end{corollary}

This result was used many times, but for the first time, apparently, it was obtained in \cite{Chemin}.

\medskip

Now we apply Lemma \ref{L1} to the case  \eqref{phi}, where the integration is over ${\mathcal B}_R$, the support of $\phi$, and all integrals converge.

Let us denote $\Omega_1=\left(x\Bigl||{\bf x}|\le
R-\frac{R}{k}\right),\,$ $\Omega_2=\left(x\Bigl| R-\frac{R}{k}<|{\bf
x}|\le R\right).$ Thus, ${\mathcal B}_R=\Omega_1\cup\Omega_2.$

\begin{corollary}
 For $\Phi$ given as \eqref{phi}
\begin{eqnarray}\nonumber
G'_\phi(t)=
 -2B\int\limits_{\Omega_1}
({\bf V}\cdot{\bf x})\,\rho\,dx+ 2C\int\limits_{\Omega_2}\frac{|{\bf
x}|-R}{|{\bf x}|} ({\bf V}\cdot{\bf x})\,\rho\,dx,\label{(1.6.2)}
\end{eqnarray}

\begin{eqnarray}
G''_\phi(t)=-2B\left(\int\limits_{\Omega_1}\rho |{\bf
V}|^2\,dx+n\int\limits_{\Omega_1}
p\,dx\right)+\nonumber\end{eqnarray}
\begin{eqnarray}+2C\left(\int\limits_{\Omega_2}\rho |{\bf
V}|^2\,dx+n\int\limits_{\Omega_2} p\,dx-R\int\limits_{\Omega_2}\rho
\frac{|\sigma|^2}{|{\bf
x}|^3}\,dx-(n-1) \, R \,\int\limits_{\Omega_2}\frac{p}{|{\bf
x}|}\,dx\right).\label{(1.6.3a)} \end{eqnarray}
\end{corollary}

\medskip

\subsection{Estimates}

Let us estimate $G''_\phi(t)$ from below. As follows from (\ref{(1.6.3a)}),
$$G''_\phi\ge
-2B\left(\int\limits_{\Omega_1}\rho |{\bf
V}|^2\,dx+n\int\limits_{\Omega_1} p\,dx\right)$$
$$+2C\left(\left(1-\frac{R}{R-\frac{R}{k}}\right)
\int\limits_{\Omega_2}\rho |{\bf
V}|^2\,dx+\left(n-(n-1)\frac{R}{R-\frac{R}{k}}\right)\int\limits_{\Omega_2}
p\,dx\right)$$
$$=2C\frac{k-n}{k-1}\int\limits_{B_R}
p\,dx-2B\left(\int\limits_{\Omega_1}\rho |{\bf
V}|^2\,dx+n\int\limits_{\Omega_1} p\,dx\right)$$
$$-2C\left(\frac{k-n}{k-1}\int\limits_{\Omega_1}
p\,dx+\frac{1}{k-1}\int\limits_{\Omega_2}\rho |{\bf
V}|^2\,dx\right).$$ Only the first term in the last expression is positive.
In order to estimate it from below, we note that from H\"older's inequality, the entropy  equation 
\eqref{e3} and the equation of state \eqref {e4} imply
\begin{equation}\label{Ggam}(G_\phi)^\gamma(t)\le \int\limits_{{\mathcal B}_R}\rho^\gamma\,dx\,
\left(\int\limits_{{\mathcal B}_R}\phi^{\frac{\gamma}{\gamma-1}}\,dx\right)^{\gamma-1}\le \exp\{-\inf\limits_{{\mathcal B}_R} S_0\}\int\limits_{{\mathcal B}_R}p\,dx\,
\left(\int\limits_{{\mathcal B}_R}\phi^{\frac{\gamma}{\gamma-1}}\,dx\right)^{\gamma-1}.
\end{equation}
 With
$A_1=(k-n)\,\exp\{\inf\limits_{{\mathcal B}_R}
S_0\}\,\left(\int\limits_{{\mathcal B}_R}\phi^{\frac{\gamma}{\gamma-1}}(t)\,dx\right)^{1-\gamma},$
taking into account the connection between the constants $ B $ and $ C $ and the expression for the total energy, we finally get
\begin{proposition} \label{p1} For case I
\begin{equation}
G''_\phi(t)\ge 2B\left(A_1G_\phi^\gamma(t)-A_2\mathcal
E\right),\label{(1.6.5)}
\end{equation}
 with
$A_2=\max\{2,(\gamma-1)k\}.$ For case II, $R\le {\mathcal R}_0$, the total energy $\mathcal
E$ should be replaced by $\mathcal
E_{B_{{\mathcal R}(t)}}$.
\end{proposition}

\bigskip

\begin{proposition} \label{p2} For case I
\begin{equation}\label{Ep} E_p(t)\ge \frac{A_1}{(\gamma-1)(k-n)} G^\gamma_\phi(t),
\end{equation}
 \begin{equation}\label{Ggam1} G_\phi(t)\le G^+:=
\left((\gamma-1)(k-n)\frac{\mathcal E}{A_1}\right)^\frac{1}{\gamma}.
\end{equation}
For case II, $R\le {\mathcal R}_0$,  the potential energy $
E_p$ should be replaced by $
{E_p}_{{\mathcal B}_{{\mathcal R}(t)}}$.
\end{proposition}
\proof
Estimates \eqref{Ep} and \eqref{Ggam1} follow from  \eqref{Ggam} directly.

\medskip

\bigskip

Further we need an auxiliary result.
\begin{proposition} \label{L2}
 Let $y_-(t)$ and $y_+(t)$ be classical solutions of equation $$y_-^{''}=K_1 y_-^\gamma-K_2(t)$$ and inequality
 $$y_+^{''}\ge K_1 y_+^\gamma-K_2(t),$$ respectively, where
  $\gamma> 1$, $K_1\ge 0$ are constants and $K_2(t)$ is a smooth function. If the initial data are such that $0<y_-(0)\le y_+(0)$,  $y'_-(0)\le y'_+(0)$, and  $y_-(t)$ keeps positivity for  $t\in (0,T]$,
  then $$y_-(t)\le y_+(t)$$
  for  $t\in [0,T]$.

\end{proposition}
\proof
Let $z=y_+-y_-$. Then $z$ solves inequality $$z^{''}\ge K_1 (y_+^\gamma-y_-^\gamma),$$
and $z(0)\ge 0,$ $z'(0)\ge 0.$ According to the Bernoulli inequality we have
 $$z^{''}\ge K_1 (((y_+-y_-)+y_-)^\gamma-y_-^\gamma)=K_1 y_-^\gamma \left(\left(\frac{z}{ y_-}+1\right)^\gamma-1\right)\ge \gamma K_1 y_-^{\gamma-1}z,$$
 what implies $z\ge 0$ as long as  $y_->0$. $\Box$

 \medskip

\begin{proposition} \label{p4} For case I on classical solutions of \eqref{e1}--\eqref{e4} the following inequality holds:
\begin{equation*}\label{Gprim}
(G'_\phi(t))^2\le 8 C G_\phi(t)\left(\mathcal E - \frac{A_1}{(\gamma-1)(k-n)} G^\gamma_\phi(t)\right).
\end{equation*}
For case II, $R\le {\mathcal R}_0$, the total energy $\mathcal
E$ should be replaced by $\mathcal
E_{{\mathcal B}_{{\mathcal R}(t)}}$.
\end{proposition}
\proof
Equality \eqref{Gprphig} implies
\begin{eqnarray}\label{Gprphi}
 (G'_\phi(t))^2=\left(\int\limits_{{\mathbb R}^n}\frac{\phi'(|{\bf x}|)}{|{\bf x}|}
({\bf V}\cdot{\bf x})\,\rho\,dx\right)^2=
\left(\int\limits_{{\mathbb R}^n}\frac{\phi'(|{\bf x}|)}{|{\bf x}|}
({\bf V}\cdot{\bf x})\,\rho^{\frac{1}{2}}\,\rho^{\frac{1}{2}}\,\phi^{\frac{1}{2}}\,\phi^{-\frac{1}{2}}\,dx\right)^2\nonumber\\
\le G_\phi(t) \int\limits_{{\mathbb R}^n}\frac{(\phi'(|{\bf x}|))^2}{|{\bf x}|^2\,\phi(|{\bf x}|)}
|({\bf V}\cdot{\bf x})|^2\rho\,dx\le 2 G_\phi(t)\, E_k(t)\,\max\limits_{x\in {\mathcal B}_R}  \left|\frac{(\phi'(|{\bf x}|))^2}{\phi(|{\bf x}|)}\right|=\nonumber\\
 8 C G_\phi(t)\, E_k(t)\,= 8 C G_\phi(t)\, (\mathcal E - E_p(t))\, \le \left(\mathcal E - \frac{A_1}{(\gamma-1)(k-n)} G^\gamma_\phi(t)\right).\nonumber
\end{eqnarray}
Here we take into account that $\max\limits_{{\mathbb B}_R}  \left|\frac{(\phi'(|{\bf x}|))^2}{\phi(|{\bf x}|)}\right|=\frac{4k}{R^2}=4C.$
At the last step we use estimate \eqref{Ep}.

 \section{ Proof of  Theorem \ref{T1}}\label{S3}

 Let us assume that the solution is smooth and denote $$M_1(t)=\sup\limits_{x\in \Omega_2}(\rho|{\bf V}|^2+n p),\quad
  M_2(t)=\sup\limits_{|x|=R}\int\limits_0^t\left( \rho \frac{|{\bf V}|^2}{2}+\frac{\gamma}{\gamma-1}p\right)|{\bf V}| d\tau, $$
  and $${\mathcal E}_{\Omega_1}={ E}_{k,\Omega_1}+{\mathcal E}_{p,\Omega_1}=\int\limits_{\Omega_1}\left(\frac{\rho |{\bf
V}|^2}{2}+\frac{1}{\gamma-1}\,p\right)\, dx.$$  Then \eqref{(1.6.3a)} implies
\begin{eqnarray}
G''_\phi(t)\le -2B\left(\int\limits_{\Omega_1}\rho |{\bf
V}|^2\,dx+n\int\limits_{\Omega_1}
p\,dx\right)+ 2C M_1(t) |\Omega_2|
\nonumber\\
=-2B\left(\delta {\mathcal E}_{\Omega_1}+\big|n(\gamma-1)-2\big| \, \tilde E - \omega_n \,M_1(t)\,(k-1) R^n \left(1- \frac{(k-1)^n}{k^n}\right)\right),\label{Gpp}
\end{eqnarray}
where $\tilde E= { E}_{k,\Omega_1},$ for $\gamma\le 1+\frac{2}{n}$ and  $\tilde E= { E}_{p,\Omega_1},$ otherwise.

Since $ |\Omega_2| \to 0$, as $k\to \infty$, then $M_1(t)\to \sup\limits_{|x|=R}(\rho|{\bf V}|^2+n p),$ as $k\to \infty$.

Further,
\begin{equation*}\label{DE}
 \frac{d {\mathcal E}_{\Omega_1}}{dt} =-\int\limits_{\partial \Omega_1}\left( \rho \frac{|{\bf V}|^2}{2}+\frac{\gamma}{\gamma-1}p\right) ({\bf V}\cdot\nu) dS,
\end{equation*}
$\nu$ is a unit normal vector to $\partial \Omega_1$, therefore
\begin{equation}\label{DE}
 {\mathcal E}_{\Omega_1}(t) \ge {\mathcal E}_{\Omega_1}(0)- M_2(t) |S_R|,\quad |S_R|=n\omega_n R^{n-1}.
\end{equation}

Thus, \eqref{Gpp} and \eqref{DE} imply
\begin{eqnarray}\nonumber
G''_\phi(t)\le
-2B\left(\delta {\mathcal E}_{\Omega_1}(0) - \omega_n R^n \,\frac{k-1}{{k^n}} (k^n- (k-1)^n) \, M_1(t)-
\delta n\omega_n R^{n-1} M_2(t)
\right).\label{Gpp1}
\end{eqnarray}

Since $ {\mathcal E}_{\Omega_1}(0)\to {\mathcal E}$ as $R\to\infty$, and  $(k-1) \left(1- \frac{(k-1)^n}{k^n}\right)\to n,$ as $k\to \infty$, then
condition \eqref{condT0} implies that for sufficiently large $k$ and $R$ there exists a positive constant $c $ such that
 $$ \delta {\mathcal E}_{\Omega_1}(0) - \omega_n R^n \,(k-1) \left(1- \frac{(k-1)^n}{k^n}\right) \, M_1(t)-
\delta n\omega_n R^{n-1} M_2(t)
>c$$
  and
$G''_\phi(t)\le -2B c.$

The last inequality contradicts  the non-negativity of $G_\phi(t)$ for any choice of initial data. Theorem \eqref{T1} is proved. $\square$

\medskip
{\it Remarks.}

1. There are many globally smooth solutions with finite mass and energy, for the methods of their construction see  \cite {Rozanova_HYP}, \cite {Rozanova_class}. Below we  investigate what prevents condition \eqref {condT0}  from being met for these solutions.

We dwell on the simplest case   $n=1$, $\gamma=3$, $\delta=2=(\gamma-1) n$ and
 notice that \eqref{e2} -- \eqref{e4} imply
 \begin{equation}\label{pp}
   \partial_t p +
   {\bf V}\cdot\nabla_{\bf x} p+
    \gamma\, p \, {\rm div}_{\bf x} {\bf
V}  =  0.
 \end{equation}

We select the initial data as
\begin{eqnarray}\label{solCD}
{\bf V}_0= a_0{\bf x},\quad \rho_0=
\frac{1}{(1+|{\bf
x}|^2)^{2}},\quad p_0=\frac{1}{(1+|{\bf
x}|^2)
}.
 \end{eqnarray}
Let as denote $\psi(t)= \sqrt{(2+a^2_0) t^2 + 2 a_0 t + 1}$. It can be readily checked that
\begin{eqnarray}\label{sol}
  {\bf V}=\frac{\psi'(t)}{ \psi(t)}{\bf x},\quad \rho =\frac{\psi^{3}(t)}{(\psi^2(t)+x^2)^2},\quad
  p =\frac{1}{\psi(t)\,(\psi^2(t)+x^2)}
  \end{eqnarray}
is a solution to the corresponding Cauchy problem \eqref{e1}, \eqref{e2},  \eqref{pp}, \eqref{solCD}.

For this solution the moment of mass $G(t)=\frac12 \int\limits_{\mathbb R} \rho |{\bf x}|^2 dx=\frac{\pi}{4}\psi^2(t)$ exists,
the total mass is
$m=\frac{\pi}{2}$, the total energy ${\mathcal E}= \frac{\pi}{4}(2+a_0^2)$, $E_p(t)=\frac{\pi }{2\psi^2(t)}\to 0$, $E_k(t)=\frac{\pi (\psi'(t))^2}{\psi(t)}\to {\mathcal E}$, as $t\to \infty$.

Here
\begin{equation*}\label{condT0}
\lim_{R\to\infty}   \sup\limits_{{|x|=R}}\left[\int_0^t\left(\frac12 \rho|{\bf V}|^2+ \frac{\gamma}{\gamma-1} p\right)\,|{\bf V}|\,dt \right]= {{\mathcal E}}= \frac{\pi}{4}(2+a_0^2)
\end{equation*}
\begin{equation*}
\lim_{R\to\infty}  {R} \sup\limits_{{|x|=R}}\left[\rho ({\bf V}|^2+  p \right]= \frac{3\sqrt{3}}{8}(2+a_0^2),
\end{equation*}
therefore \eqref{condT0} is not satisfied. In fact, the reason for this is that although for any fixed $ t $ the energy $E_R(t)\to \mathcal E$, as $R\to \infty$, $E_R(t)\to 0$,  as $t\to \infty$, for arbitrary large fixed $R$.

\medskip

2. In \cite{Cho2006} the following sufficient condition for the blow-up for solutions with a finite moment of mass $G(t)=\frac12 \int\limits_{\mathbb R} \rho |{\bf x}|^2 dx$ was obtained (valid both for compressible Euler and Navier-Stokes
equations):
\begin{equation}\label{Cho}
  \limsup\limits_{t\to\infty} \left\|t \frac{({\bf V}\cdot {\bf x})}{(1+|{\bf x}|^2)} \right\|<1.
\end{equation}
Let us notice that for solution \eqref{sol} the limit in \eqref{Cho} is exactly 1.

Condition \eqref{condT0} is different in its nature: it controls the behaviour of solution as $|x|\to \infty $ and do not prescribes any additional requirement on $\rho$. In fact, \eqref{condT0} can be applied even to solutions with infinite mass.

 \section{ Proof of  Theorem \ref{T2}}\label{S4}

\subsection{Lower estimates of $ G_\phi$}

 Let $z(t)$ be the solution to problem \eqref{(1.6.5e)}.
Denote $z'=q(z)$, then  we get
\begin{eqnarray}\nonumber
q^2= 4B\left(\frac{A_1}{\gamma+1}  z^{\gamma+1}-A_2\mathcal
E z\right)+c,\label{(1.6.5e1)}\\  c=( G'_\phi(0))^2 -4B\left(\frac{A_1}{\gamma+1}   G_\phi^{\gamma+1}(0)-A_2\mathcal
E G_\phi(0)\right). \label{KK}\nonumber
\end{eqnarray}

We introduce function
\begin{equation}\label{cond}
f(z)= 4B\left(\frac{A_1}{\gamma+1}  z^{\gamma+1}-A_2\mathcal
E z\right)+c,
\end{equation}
and notice that
 $q=\pm \sqrt{f(z)}$ is the phase curve of \eqref{(1.6.5e)} (the phase portrait is shown in Figure 1).

The formal quadrature is
\begin{equation}\label{solz}
  t=\pm \int \frac{1}{\sqrt{f(z)}} \, dz.
\end{equation}

The function $f(z)$ is positive at the point $ G_\phi(0)$, it tends to $+\infty$ as $z\to+\infty$, and can have
two, one or no roots on $( G_\phi(0),\infty)$. The only  minimum  point is at $z=(\frac{A_2\mathcal
E}{A_1})^{\frac{1}{\gamma}}:=G^{++}.$ Let us show that
\begin{equation}
G_\phi^\gamma(0)<
G^{++},\label{ie1}
\end{equation}
that is $G^{++}\in ( G_\phi(0),\infty).$ Indeed,
 the right hand side of  \eqref{ie1} can be estimated from below as
\begin{eqnarray}
\frac{A_2}{A_1}\mathcal E\ge \frac{A_2}{A_1}E_p(0)=
\frac{A_2}{A_1(\gamma-1)}\int\limits_{{\mathcal B}_R}p\,dx>\nonumber\\
 \frac
{\max\{2,(\gamma-1)k\}}{(k-n)(\gamma-1)}\,\exp\{\inf\limits_{{\mathcal B}_R}
S_0\}\,
\left(\int\limits_{{\mathcal B}_R}\phi^{\frac{\gamma}{\gamma-1}})\,dx\right)^{\gamma-1}\int\limits_{{\mathcal B}_R}\rho^\gamma\,dx
.\label{(1.6.15)} \end{eqnarray}

The left hand side of \eqref{ie1}  can be estimated from above as \eqref{Ggam}. Therefore, assuming inequality, opposite to  \eqref{ie1},
from
\eqref{(1.6.15)} and \eqref{Ggam} we get
$(\gamma-1)(k-n)>\max\{2,(\gamma-1)k\},$ it is a contradiction.

This contradiction, in particular, implies
\begin{equation}\label{GGG}
G_\phi(t)<G^+< G^{++},
\end{equation}
$G^+$ is defined in \eqref{Ggam1},
for all $t$, since we can shift the start time to a new point.
We distinguish two cases:
\begin{enumerate}
  \item
  $f(z)$ has no roots for $z\in ( G_\phi(0),\infty)$. It means that $z$ increases unboundedly;
  \item
  $f(z)$ has two roots (or one root of multiplicity 2)  for $z\in ( G_\phi(0),\infty)$. We denote $z_*$
  the smaller positive root.
  \end{enumerate}

Proposition \ref{L2} means that if $G_\phi(t)$ is  a solution to differential inequality \eqref{(1.6.5)}, then
 $ G_\phi(t)\ge z(t),$ for $t\in [0,T]$ such that $z(t)>0$.

 In case 1, when $f(z)>0$ for all $z>G_\phi(0)$,
 \eqref{solz} implies $t\to \infty $ as $z\to \infty$, therefore $T=\infty.$
 In case 2 the time $T$ can be computed as
 \begin{eqnarray}\label{time1}
 T&=&- \int\limits_{ G_\phi(0)}^{0} \frac{1}{\sqrt{f(g)}} \, dg,\quad\mbox{if}\quad G_\phi'(0)\le 0, \quad\mbox{and}\\
  T&=&\int\limits_{ G_\phi(0)}^{z_*} \frac{1}{\sqrt{f(g)}} \, dg+ \int\limits_{0}^{ z_*} \frac{1}{\sqrt{f(g)}} \, dg,\quad\mbox{if}\quad
 G_\phi'(0)> 0,\label{time2}\end{eqnarray}
see Figure 1 for the configuration of the phase curves.
 Thus, for any fixed
 $\bar t\in [0,T]$ we have
 $$
 z(\bar t)\le  G_\phi(\bar t) \le \sup\limits_{x\in B_R} \rho(\bar t,x) \int\limits_{{\mathcal B}_R}\phi(|{\bf x}|)\, dx= \sup\limits_{x\in B_R} \rho(\bar t,x)\, K.
 $$
 Since $z(t)=z_-(t)$, we obtain \eqref{T21}. Inequality, opposite to \eqref{T21}, leads to a contradiction and implies that the solution
 loses smoothness before  the time $\bar t.$

\subsection{Upper estimates of $ G_\phi$}
Evidently,
\begin{equation*}
0<G_\phi(t)< m.\label{(1.6.4)} \end{equation*}
 If we denote
\begin{eqnarray}\label{cond00}
z_+:=\min\{ m, G^+   \},
\end{eqnarray}
where  $G^+$ is defined in \eqref{Ggam1}, then from \eqref{GGG} we have  $ G_\phi(t)<z_+$ for all $t\ge 0.$

Thus, if the solution keeps smoothness,  for any fixed
 $\bar t \ge 0 $ we have
 $$
 z_+\ge  G_\phi(\bar t) \ge \inf\limits_{x\in B_R} \rho(\bar t,x) \int\limits_{{\mathcal B}_R}\phi(|{\bf x}|)\, dx= \inf\limits_{x\in B_R} \rho(\bar t,x)\, K.
 $$
 it implies \eqref{T22}. Inequality, opposite to \eqref{T22}, leads to a contradiction and implies that the solution
 loses smoothness before  the time $\bar t.$

 Thus, Theorem \ref{T1} is proved.
 $\Box$

\subsection{"Phantom" condition for singularities formation}

It is easy to find conditions leading to a contradiction with  the upper and lower bounds of  $ G_\phi(t)$.
Indeed, if we are in the case 1, when $f(z)$ has no roots as $z\in ( G_\phi(0),\infty)$, then  $z$ increases unboundedly and reaches $z_+$ in a finite time. If we are in the case 2  and the lower positive root $z_*$ of $f(z)$ on the semi-axis $( G_\phi(0),\infty)$ is greater than $z_+$ ($f'(z_*)\le 0$), then $z$ also reaches $z_+$ in a finite time. In the both cases it implies a contradiction with inequality
$ G_\phi(t)<z_+$.

Thus, the finite-time singularity condition in terms of the initial data looks like
 \begin{equation}
 G'_\phi(0)>0, \quad f(z_+)\ge 0,\label{ie2}
\end{equation}
with $z_+$ defined in \eqref{cond00}.

However, we could neither find the initial data satisfying \eqref{ie2}, nor prove the inconsistency of this condition.

\begin{figure}[h!]
\centerline{\includegraphics[width=0.5\columnwidth]{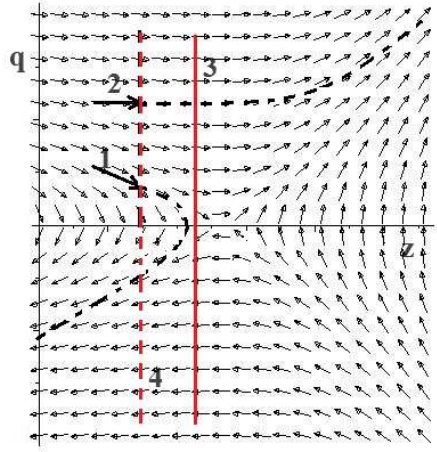}}
\caption{ The phase portrait of \eqref{(1.6.5e)}. Trajectory corresponding  to  the data such the  $f(z)$ has two roots, realistic (curve 1) and no roots, hypothetic (curve 2) on $(G_\phi(0), \infty)$. Upper bounds $z_+$ for the realistic data (line 3) and  the hypothetic data $f(z_+)>0$, leading to  singularity formation  (line 4).}
\end{figure}

 \subsection {Dynamics of momentum $G_\phi$}
 As follows from Proposition \ref{p4}, all possible values of $(z,q)$, where $z=G_\phi, q=G'_\phi$, are bounded by the graph of function
 \begin{equation}\label{ff}
q^2=8 C z\left(\mathcal E - \frac{A_1}{(\gamma-1)(k-n)} z^\gamma\right).
\end{equation}
Moreover, for any fixed $\bar t\in [0,T]$ (see \eqref{time1}, \eqref{time2}) all possible values of $G\phi(\bar t)$ lie to the right of the respective value of $z(\bar t)$ on the graph
 \begin{eqnarray*}\label{eq}
 q^2=f(z)=4B\left(\frac{A_1}{\gamma+1}  (z^{\gamma+1}-G_\phi^{\gamma+1}(0)) -A_2\mathcal
E(z-G_\phi(0))\right)+(G'_\phi(0))^2.
 \end{eqnarray*}

Possible dynamics of momentum $G_\phi$ provided the solution keeps smoothness is presented in Figure 2 on example of globally smooth solution \eqref{sol}.

\begin{figure}[h!]
\centerline{\includegraphics[width=0.5\columnwidth]{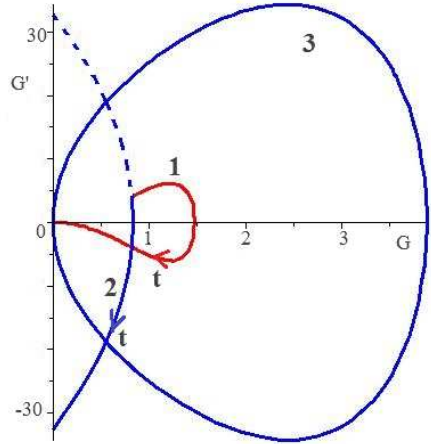}}
\caption{Dynamics of $G_\phi(t)$ on the plane $(G_\phi, G_\phi')$ on example of smooth solution \eqref{sol} ($n=1$, $\gamma=3$, $R=1$, $k=2$, $a_0=-7$) (curve 1); for the initial data \eqref{solCD}: the position of  lower bounds of $G_\phi$ at a fixed $t$ (curve 2, given as \eqref{(1.6.5e)}) and the limits  of domain in which the points $(G_\phi, G_\phi')$ lie (curve 3, given as \eqref{ff}).}
\end{figure}

\bigskip
\section{ Proof of  Theorem \ref{T3} }\label{S5}

First of all, we note that the assumptions of Theorem \ref{T3} imply $T_2>0$ and ${\mathcal B}_{{\mathcal R}(0)}\subset \Omega_1=
{\mathcal B}_{R-\frac{R}{k}}$.

Further, from \eqref{(1.6.3a)} we get
\begin{eqnarray}
G''_\phi(t)&= & -2B\left(\int\limits_{\Omega_1}\rho |{\bf
V}|^2\,dx+n\int\limits_{\Omega_1}
p\,dx\right)+ 2C \bar p \left(\int\limits_{\Omega_2}\,dx - (n-1)R \int\limits_{\Omega_2}  \frac{1}{|{\bf x}|}\,dx      \right)
\nonumber\\
&=&-2B\left(\delta {\mathcal E}_{\Omega_1}+\big|n(\gamma-1)-2\big| \tilde E - \bar p \omega_n n R^n \left(\frac{k-1}{k}\right)^n\right)\nonumber\\
&\le & - 2B \left( \delta e(t) +\frac{\delta}{\gamma-1} \bar p  \int\limits_{\Omega_1}\,dx  - \bar p \omega_n n R^n \left(\frac{k-1}{k}\right)^n  \right)\nonumber\\
&=&- 2B \left( \delta e(t) +
  \omega_n  R^n \bar p  \left(\frac{k-1}{k}\right)^n
\left(\frac{\delta}{\gamma-1}-n\right)  \right)\nonumber\\&=&- 2BN,\label{GppT3}
\end{eqnarray}
where we denoted, as in the proof of Theorem \ref{T1}, $\tilde E= { E}_{k,\Omega_1},$ for $\gamma\le 1+\frac{2}{n}$ and  $\tilde E= { E}_{p,\Omega_1},$ otherwise.

If $N>0$, then \eqref{GppT3} signifies that there exists a positive $T=T_1$ such that $G_\phi(T)\le 0$.
If $N<0$, then $T_1$ is the smallest positive root of polynomial  $G_\phi(t)=-\delta B N t^2 + G_\phi' (0) t +G_\phi (0)$ in case such a root exists ($G_\phi' (0)$ should be negative). If $N=0$, then a limit pass as $N\to 0$  in $T_1$ gives $T_1=-\frac{G_\phi (0)}{G_\phi' (0)}$. In all cases, we obtain a contradiction with the positivity of $G_\phi(t)$.

The condition \eqref {T3cond} means that $ G_\phi (t) $ becomes non-positive inside $ \Omega_1 $, and the perturbation does not have time to reach its boundary. Thus, Theorem \ref{T3} is proved. $\Box$

\medskip

We can see that if $G_\phi(t) $ has positive roots, then, by choosing $\sigma$ small enough, we can always obtain the implementation of \eqref{T3cond}.


\section{ Proof of  Theorem \ref{T4}}\label{S6}

Let us note that $Q_\phi(t)=G_\phi(t)-\bar\rho K.$

First of all, we obtain the upper bound of $Q_\phi$ as
\begin{eqnarray}\nonumber
Q_\phi(t)&=&m(0)+\int\limits_{{\mathbb R}^n}(\rho-\bar \rho) (\phi(|{\bf
x}|)-1)\,dx \\&=&m(0)+ \int\limits_{{\mathbb R}^n}\rho (\phi(|{\bf
x}|)-1)\,dx +\bar\rho\int\limits_{{\mathbb R}^n} (1- \phi(|{\bf
x}|))\,dx \nonumber \\&\le &m(0)-\bar\rho K+\bar\rho \omega_n ({\mathcal R}_0+\sigma t)^n=Q^+(t),\label{QB}\nonumber
\end{eqnarray}
therefore $Q_\phi$ grows at most polynomially.

Let us notice that $G''_\phi=Q''_\phi$. Since  $\frac{Q_\phi}{\bar\rho K}>-1$, taking into account the Bernoulli inequality we get
\begin{eqnarray}\nonumber
(G_\phi)^\gamma&=&(Q_\phi+\bar\rho K)^\gamma=(\bar\rho K)^\gamma\left(\frac{Q_\phi}{\bar\rho K}+1\right)^\gamma\\&\ge& (\bar\rho K)^\gamma\left(1+\gamma \frac{Q_\phi}{\bar\rho K}\right)=
 (\bar\rho K)^{\gamma}+\gamma (\bar\rho K)^{\gamma-1} Q_\phi.\nonumber
\end{eqnarray}
Thus, Proposition \ref{p1} implies
\begin{equation}\nonumber
Q''_\phi(t)\ge \varkappa^2 Q_\phi(t)     +P_n(t),
\label{QQ}
\end{equation}
where $\varkappa^2 = 2\gamma B A_1 (\bar\rho K)^{\gamma-1}$ is a positive constant, $P_n(t)= 2 B (A_1 (\bar\rho K)^{\gamma} -A_2 (e(0)+\bar p \,\omega_n \,({\mathcal R}_0+\sigma t)^n))$ is a polynomial of order $n$ with respect to $t$.
Further, if $z(t)$ is a solution to the Cauchy problem
\begin{equation}
 z''(t)= \varkappa^2  z(t)     +P_n(t),\quad  z(0)= Q_\phi(0), \quad z'(0)=Q_\phi'(0),
\label{QQCP}
\end{equation}
then by the comparison theorem we have $Q_\phi(t)\ge z(t):= Q_-(t)$ for all $t>0$. The solution to \eqref{QQCP} can be easily found:
\begin{equation*}
z_\phi= C_1 e^{\varkappa t} + C_2 e^{-\varkappa t} +\tilde P_n(t),
\label{Qbar}
\end{equation*}
$\tilde P_n(t)$ is a polynomial of order $n$ with respect to $t$, $C_1$ and $C_2$ depend on $Q_\phi(0)=G_\phi(0)-\bar\rho K,$ $Q_\phi'(0)=G_\phi'(0)$.

Thus, as long as the solution remains smooth, we have $ Q_-(t)\ge Q(t)\ge Q^+(t).$
Proceeding similarly to the proof of  Theorem \ref{T2}, we obtain
 \eqref{T41} and \eqref{T42}.

Thus, Theorem \ref{T4} is proved.
$\Box$

\bigskip

\section{Discussion}\label{S7}

1.
We see that the sufficient conditions for the formation of a singularity in a finite time, indicated in Theorems \ref{T1} and \ref{T3}, use some requirements for the behavior of the solution as $|{\bf x}|\to \infty$. In case I, this is a limitation on too much growth of velocity,  in case II, this is a limitation for too fast propagation of the support of perturbation.
As far as we know, all sufficient conditions of this kind contain one or another restriction. A natural question arises: is it possible to formulate integral sufficient conditions only in terms of initial conditions, without additional restrictions?

From the proof of Theorem \ref{T2} we see that
 if  there exist data  satisfying  \eqref{ie2},
then the solution of \eqref{e1} -- \eqref{e4}, \eqref{cd}  with  total mass $m$ and total energy $\mathcal E$
cannot be smooth for all $t>0$.

However, we have not been able to find the initial data satisfying \eqref {ie2}, and it is not obvious that they are inconsistent, i.e. that the class of the corresponding initial data is empty.

In the proof of Theorem \ref{T4} one can also find such "empty" sufficient conditions, for example, the existence of $T$ such that $Q^-(T)>Q^+(T)$.

\medskip

2. In the statement of Theorem \ref{T2}   we can consider another lower bound $z_-(t)$, similarly to the lower bound used in Theorem \ref{T4} (given by the solution of linear equation \eqref{QQCP}). In this case, we have no restriction on $T$ from above. Nevertheless, this estimate is rougher, moreover, starting from some $ t $ this new lower bound $ z_- (t) $ becomes negative, and therefore this estimate is trivial (it does not make sense for a positive $ G_\phi (t) $). In Theorem \ref {T4}, it would seem that such problems do not arise, since $ Q (t) $ can be negative. However, if we examine the growth of $Q(t)$ (polynomial) and $Q_-(t)$ (exponential), we will see that starting from some sufficiently large $t$ the lower estimate of $Q(t)$  still becomes trivial.

\medskip

3. It would seem that all these theorems can be applied to  the compressible Navier-Stokes equations.
 However, this is not the case. The point is that in order for the integrals over the gluing sets
 $ | {\bf x} | = R $ and $ | {\bf x} | = R- \frac{R} {k} $ the $ C^1 $ - smoothness of the function $ \phi ({\bf x})$ is not enough to vanish .

\medskip
4. Theorem  \ref{T3} can be applied under the same conditions that the original theorem of    \cite{Sideris}, therefore it is interesting to compare the results. Analysis shows that Theorem \ref {T3} covers a broader situation. Indeed, in \cite{Sideris} the condition
$$\int\limits_{{\mathbb R}^n} ( {\bf V}_0 \cdot {\bf x})\,\rho_0 \, dx>0 $$ is the main requirement, therefore, this result concerns initially highly divergent flows, while the conditions of Theorem  \ref{T3} are more flexible and cover a wide variety of initial data.
Note that the method of \cite{Sideris} can be quite naturally extended to the class of initially converging flows (similarly to \cite{RozanovaJMS}, where a more complicated problem is considered for the Euler equations containing the Coriolis force).

\medskip

5. The estimates for a smooth solution inside a given region of space, made in Theorems \ref {T2} and \ref {T4}, can be useful in numerical calculations. Indeed, it is often easy to control maxima and minima of solution, while it is difficult to decide whether we are seeing a shock wave or simply a zone of high gradients. Experiments with real data show that such estimates make sense for sufficiently small $ R $, otherwise they are very rough.

\end{document}